\documentclass[11pt]{article}

\pdfoutput=1

\usepackage{amsmath}
\usepackage{amssymb}
\usepackage{enumerate}
\usepackage{amscd}

\begin{document}

\newcommand{\LP}{\mathrel{\mathbb{L}}}
\newcommand{\RP}{\mathrel{\mathbb{R}}}
\newcommand{\HP}{\mathrel{\mathbb{H}}}
\newcommand{\DP}{\mathrel{\mathbb{D}}}
\newcommand{\JP}{\mathrel{\mathbb{J}}}

\newcommand{\eop}{\hfill$\square$}

\newcommand{\be}{\begin{enumerate}}
\newcommand{\ee}{\end{enumerate}}
\newcommand{\CL}{\mathrel{\mathcal {L}}}
\newcommand{\CR}{\mathrel{\mathcal {R}}}
\newcommand{\CH}{\mathrel{\mathcal {H}}}
\newcommand{\CJ}{\mathrel{\mathcal {J}}}

\newcommand{\ZZ}{\mathbb{Z}}
\newcommand{\NN}{\mathbb{N}}
\newcommand{\LE}{\widetilde{\CL}_E}
\newcommand{\RE}{\widetilde{\CR}_E}
\newcommand{\HE}{\widetilde{\CH}_E}

\newcommand{\BV}{{\bf V}}
\newcommand{\BG}{{\bf G}}

\newcommand{\BI}{{\bf I}}

\newcommand{\BN}{{\bf N}}
\newcommand{\BM}{{\bf M}}

\newcommand{\BSL}{{\bf SL}}

\newcommand{\End}{{\rm End}\,}
\newcommand{\Ord}{{\rm Ord}\,}

\newcommand{\darrow}{\!\downarrow}
\newcommand{\np}{^e}
\newcommand{\thetap}{\theta \np}
\newcommand{\psip}{\psi \np}
\newcommand{\otheta}{\overline{\theta}}
\newcommand{\oalpha}{\overline{\alpha}}
\newcommand{\obeta}{\overline{\beta}}

\newcommand{\CM}{{\cal M}}
\newcommand{\CPM}{{\cal PM}}

\newtheorem{propo}{PROPOSITION}[section]
\newtheorem{result}[propo]{RESULT}
\newtheorem{lemma}[propo]{LEMMA}
\newtheorem{theorem}[propo]{THEOREM}
\newtheorem{example}[propo]{EXAMPLE}
\newtheorem{corol}[propo]{COROLLARY}
\newtheorem{question}[propo]{QUESTION}

\title{Almost perfect restriction semigroups}
\author{Peter R. Jones\\
Department of Mathematics, Statistics and Computer Science\\
Marquette University\\
Milwaukee, WI 53201, USA\\
peter.jones@mu.edu}

\maketitle

\begin{abstract}
We call a restriction semigroup {\em almost perfect} if it is proper and its least monoid congruence is perfect.  We show that any such semigroup is isomorphic to a `$W$-product' $W(T,Y)$, where $T$ is a monoid, $Y$ is a semilattice and there is a homomorphism from $T$ into the inverse semigroup $TI_Y$ of isomorphisms between ideals of $Y$. Conversely, all such $W$-products are almost perfect.  Since we also show that every restriction semigroup has an easily computed cover of this type, the combination yields a `McAlister-type' theorem for all restriction semigroups.  It is one of the theses of this work that almost perfection and perfection, the analogue of this definition for restriction monoids, are the appropriate settings for such a theorem.  That these theorems do {\em not\/} reduce to a general theorem for inverse semigroups illustrates a second thesis of this work: that restriction (and, by extension, Ehresmann) semigroups have a rich theory that does not consist merely of generalizations of inverse semigroup theory.   It is then with some ambivalence that we show that all the main results of this work easily generalize to encompass {\em all\/} proper restriction semigroups.

The notation $W(T,Y)$ recognizes that it is a far-reaching generalization of a long-known similarly titled construction.  As a result, our work generalizes Szendrei's description of almost factorizable semigroups while at the same time including certain classes of free restriction semigroups in its realm.\\



\end{abstract}

\section{Introduction}\label{intro}
The study of the structure of restriction semigroups has in large part consisted of generalizations of structure theorems for inverse semigroups.  For instance, the Munn representation of inverse semigroups by isomorphisms between the principal ideals of its semilattice of idempotents generalizes naturally \cite{FGG1} to representations of restriction semigroups by similar mappings of its semilattice of projections, and this generalized representation is at the very heart of our work. 
The `inductive groupoid' approach to inverse semigroups has been extended successfully to restriction semigroups by Lawson \cite{lawson}.

Somewhat complementary to the Munn representation is the McAlister theory, whereby every inverse semigroup is an idempotent-separating image of an $E$-unitary inverse semigroup, and the latter semigroups are described  as `$P$-semigroups', relative to their semilattices, greatest group images and one further structural parameter.  This theory has been generalized with some success to restriction semigroups, with $E$-unitariness replaced by `properness'.

When moving yet further from inverse semigroups, Branco, Gomes and Gould \cite{BGG, GG} introduced the notion of $T$-properness of (one-sided) Ehresmann semigroups, with respect to a submonoid $T$.  The main thesis of our work is that (returning to the realm of restriction semigroups) this idea yields narrower notions of properness that are the appropriate ones in which to prove a `McAlister-type' theory.  That this theory specializes in the case of proper inverse semigroups to a very narrow subclass we take to be a witness to our thesis, rather than the contrary.  Our results suggest that the road taken for Ehresmann semigroups in the cited papers is indeed a natural one.

To illustrate that this is not merely a conceit, we state the main result of this paper, to illustrate its simplicity. (In fact, we prove a more general theorem, applying to all proper restriction semigroups, about which more will be said below.) Recall first that for restriction semigroups, monoids play the role that groups play for inverse semigroups and that such a semigroup is {\em proper\/} if the least monoid congruence $\sigma$ meets each class of the generalized Green relations $\LP$ and $\RP$ trivially.

We call a restriction semigroup $S$ {\em almost proper\/} if it is proper and $\sigma$ is perfect (meaning that the product of classes is again a class).  The reason for the qualifier `almost' is that we term a restriction {\em monoid\/} $M$ {\em perfect\/} if, further, each $\sigma$-class has a greatest element, that is, $M$ is also an {\em $F$-restriction monoid\/}.  (The connections with $T$-properness will be made below.)

Now suppose $T$ is a monoid, $Y$ is a semilattice and there is a representation of $T$ in the inverse semigroup $TI_Y$ of isomorphisms between ideals of $Y$. This corresponds to a partial (right) action $(t, e) \mapsto e^t$, for $t \in T$ and $e \in \Delta t$, the domain of (the image of) $t$, where $e^t \in \nabla t$, the range of (the image of) $t$.  Let 
\[W(T, Y) =   \{(t,e^t) \in T \times Y :  e \in \Delta t\}\quad  (=\{(t,f) \in T \times Y : f \in \nabla t\}).\]

The multiplication in this `$W$-product' is defined by
\( (t,e^t)(u, f^u) = (tu, (e^t f)^u)\);
unary operations are defined by
\((t,e^t)^+ = (1,e)\)  and \( (t,f)^* = (1,f).\)

\begin{theorem}\label{main} (Theorems~\ref{X semigroup} and \ref{converse}, Corollary~\ref{T-proper cover}) Every restriction semigroup has an almost perfect (projection-separating) cover.  Every semigroup $W(T,Y)$ is almost perfect and, conversely, every almost perfect restriction semigroup $S$ is isomorphic to a semigroup of that form. 

Every restriction monoid has a perfect monoidal cover.  If the monoid $T$ acts by isomorphisms between {\rm principal} ideals of a semilattice $Y = Y^1$, then $W(T,Y)$ is a perfect monoid and, conversely, every perfect restriction monoid $M$ is isomorphic to such a monoid, where $Y$ is the semilattice $P_M$ of projections of $M$, $T = M/\sigma$, and the action is induced by the (generalized) Munn representation of $M$.
\end{theorem}

We must emphasize that our construction, although using the same notation, is more general than the original construction \cite{FG, GS, S1}, which corresponds precisely to the case that the representation is by {\em fully defined\/} isomorphisms between ideals of $Y$, that is, by injective endomorphisms whose ranges are ideals of $Y$. See the discussion following Theorem~\ref{X semigroup}, and its application in Proposition~\ref{factorizable}.

The connection with the term $T$-proper mentioned earlier is a central part of this work.  We call a proper monoid $M$ {\em $T$-proper\/} if it contains a `plain' submonoid $T$ such that $M = P_M T$ and $T$ is separated by $\sigma$; equivalently, each element of $M$ can be uniquely expressed in the form $e t$, where $e \in P_M, t \in T$.  In fact, $T$-properness is equivalent to perfection (Proposition~\ref{T-proper}). For a restriction {\em semigroup\/}, almost perfection is equivalent to {\em almost $T$-properness\/} (Proposition~\ref{almost T-proper}), which is the property that the monoid $C(S)$ of permissible sets, in which $S$ embeds, is $T$-proper (that is, perfect).

In the semigroup case of the theorem above,  the representation is obtained as follows: $Y$ is again the semilattice of projections, $T $ is a monoid such that $S^1 = P_{S^1} T$, and the action is induced by the extension of the (generalized) Munn representation $S$ to the monoid $C(S)$.  The key role played by $C(S)$ was suggested by its role in the determination of proper left factorizable and almost factorizable semigroups by Gomes and Szendrei \cite{GS, S1}, which together with their use of the original construction $W(T,Y)$ provided some of the motivation for our techniques.

As well as applying to both the one- and two- sided versions of factorizable and almost factorizable semigroups (see Section~\ref{embedding theorem}), our construction applies to free restriction semigroups (and monoids) and to relatively free semigroups in certain varieties of restriction semigroups.  (See Section~\ref{examples}.) This application provided the second main impetus for the $T$-proper approach to our work.

In Section~\ref{embedding theorem} we show how our methods can be used to re-prove the recent result of Gomes, Gould and Hartmann \cite{gould et al}, kindly related to the author in a private communication, that any restriction semigroup can be embedded into a factorizable one.

In fact, the construction of the semigroups $W(T, Y)$ and the constructions used in the covering theorems -- and in fact a description of {\em all\/} such covers (Corollary~\ref{all covers}) -- are all instances of a very general construction 
$S_{T,R}$ (Theorem~\ref{newcovers}) which, despite its somewhat technical hypotheses, has a simple verification since the calculations are performed within a direct product $S \times T$.  The virtue of this approach is witnessed in the following section, when all the properties (associativity included) of the $W$-product follow. Further, the representation theorem by the $W$-product also then follows from the description of covers mentioned above.

In Sections~\ref{examples} and \ref{embedding theorem}, we apply the $W$-product construction to various representations.  As further evidence of the centrality of almost perfect restriction semigroups, we note there that whenever a variety of restriction semigroups is closed under taking proper covers, the free objects in the variety are necessarily almost perfect.  We also provide a unified description of left, right, and two-sided almost factorizable restriction semigroups (cf \cite{S1}).

In the penultimate section of the paper, a further level of generality, about which we have a certain ambivalence, is achieved by weakening the requirement that the map $\alpha: T \longrightarrow TI_Y$ (as in the statement of Theorem~\ref{main}) be a homomorphism.  We now must allow {\em subhomomorphisms\/}: mappings $\alpha$ such that only $  (a\alpha)(b\alpha) \leq (ab) \alpha$, referring to the natural partial order, is required.  By means of this generalization, {\em all\/} proper restriction semigroups are described via the $W$-product construction.  In fact, all the general theorems of the paper have straightforward extensions to the proper case by this means.

Our ambivalence stems from the main thesis of this paper: that the concept of perfection (for monoids) and almost perfection (for semigroups) is in reality the `optimum' one, rather than properness per se, with Theorem~\ref{main} and the examples in Sections~\ref{examples} and \ref{embedding theorem} our support for this thesis, in practical terms, along with the elegance of homomorphic actions by monoids, rather than subhomomorphic ones.

The final section specializes the general results of the previous section to inverse semigroups: the representation turns out to be essentially that of Petrich and Reilly \cite{PR}, though this played no role in the development of our work.  More importantly, we demonstrate that the covering theorems, in particular, do {\em not} extend in a meaningful way to inverse semigroups, producing a restriction semigroup even when starting from an inverse semigroup, in general.

Although we believe that this paper incorporates a new approach, there is a already a body of work on the structure of proper restriction semigroups.  Cornock and Gould \cite{CG} obtained a structure theorem for proper restriction semigroups in general, based more on the work of Petrich and Reilly \cite{PR}, cited above, than on McAlister's theorem {\it per se\/}, using as parameters a monoid that acts on both sides of a semilattice.  Of course their results and those in Section~\ref{proper} of our paper must be deducible from each other, but we leave that exercise to the reader, as we leave the exercise of specializing the results therein to the perfect and almost perfect cases.

As this paper was nearing submission, the author received a copy of the preprint \cite{K} by G. Kudryavtseva. Although, by and large, her paper is complementary to this one, in goals and approach, being primarily expressed in terms of actions and focusing largely on restriction {\em monoids\/}, there is significant overlap.  For instance, $F$-restriction monoids appear under the same name (unsurprisingly, given their genesis in $F$-inverse monoids).  Her `ultra-$F$ restriction' monoids apparently coincide with our perfect monoids and her construction $M(T,Y)$ is therefore inevitably equivalent to our generalized $W(T,Y)$ construction in the monoidal case.  She makes explicit the connection with the work of \cite{CG} and with \cite{FGG} in a way that we do not.

In Section~\ref{embedding theorem}, we shall make further allusion to the part of \cite{K} that develops embedding theorems related to the work of Szendrei.

\section{Preliminaries}\label{prelim}

We first introduce restriction semigroups more formally, along with their basic properties and related definitions.  A {\em left restriction\/} semigroup is a unary semigroup $(S, \cdot, ^+)$ that satisfies the `left restriction' identities
\noindent \[ \quad  x^+x = x; \quad \quad (x^+y)^+ = x^+ y^+;\quad  \quad x^+ y^+ = y^+ x^+; \quad  \quad xy^+ = (xy)^+ x.\]

\noindent A {\em right restriction\/} semigroup is a unary semigroup $(S, \cdot, ^*)$ that satisfies the `dual' identities, obtained by replacing $^+$ by $^*$ and reversing the order of each expression.  A {\em restriction\/} semigroup is then a biunary semigroup $(S, \cdot, ^+, ^*)$ that satisfies both sets of identities, along with $(x^+)^* = x^+$ and $(x^*)^+ = x^*$.

From these identities it follows that for all $x \in S$, $x^+$ is idempotent and $(x^+)^+ = x^+$. These idempotents are the {\em projections\/} of $S$; by duality these are also the idempotents $x^*$, $x \in S$.  Denote the set of projections by $P_S$.  Although, by the third identity, $P_S$ is a semilattice, this need by no means be true of $E_S$.  The following consequence of the identities is well known.

\begin{lemma}\label{basic lemma} Let $S$ be a restriction semigroup. Then $S$ satisfies $x^+ \geq (xy)^+$ and $(xy)^+ = (xy^+)^+$ and their duals, namely $y^* \geq (xy)^*$ and $(xy)^* = (x^*y)^*$.

\end{lemma}

Until quite recently, the term `weakly $E$-ample' was used for restriction semigroups, providing evidence of a succession of generalizations --- by the so-called York school --- of Fountain's `ample semigroups', which we will define below. 

A restriction monoid is a restriction semigroup with identity 1 and $1^+ = 1^* = 1$.  When adjoining an identity element to a restriction semigroup, in the usual way, the latter relations are assumed.  In the standard terminology, restriction semigroups $S$ with $|P_S| = 1$ are termed {\em reduced\/}. Since, in essence, they are just monoids, regarded as restriction semigroups by setting $a^+ = a^* = 1$ for all $a$, we will generally use the latter term, except in case of possible ambiguity. 

The relevant generalized Green relations may be defined as follows. In any restriction semigroup, $\RP \, =  \{(a,b): a^+ = b^+\}$, $\LP \, = \{(a,b): a^* = b^*\}$ and $\HP \, = \, \LP \cap \RP$.  It follows easily from Lemma~\ref{basic lemma} that each contains the corresponding usual Green relations, that $\RP$ is a left congruence and that $\LP$ is a right congruence. In the standard literature, these relations are usually denoted $\RE$, $\LE$ and $\HE$, respectively (more recently, e.g. \cite{S1}, without the subscripts).  The semigroup is {\em left ample\/} if $\RP$ coincides with the relation $\CR^* = \{(a,b): xa = xb \text{ if and only if } xb = yb, \text{ for all } x, y \in S^1\}$; {\em right ample\/} if the dual statement holds, and {\em ample\/} if both hold.

The {\em natural partial order\/} on $S$ is defined by $a \leq b$ if $a = e b$ for some $e \in P_S$; equivalently if $a = a^+b$.  It is self-dual, compatible with the operations on $S$ and extends the usual order on $P_S$.  Put $a \darrow = \{ b \in S : b \leq a\}$, the principal order ideal generated by $a$.  An {\em order ideal\/} of $S$ is then a nonempty subset that is closed under $\darrow$. The order ideals of a semilattice are just its ideals.

In general, the term `homomorphism' will be used appropriate to context: that is, it should respect the binary operation for `plain' semigroups, and either or both unary operations for one-sided or two-sided restriction semigroups. In the case of monoids, we shall use the qualifier `monoidal' to indicate that it should respect the identity elements. When considering subsemigroups (or submonoids), we shall usually use the qualifier `biunary' or `restriction' in the latter case when the situation might not be clear.  A biunary subsemigroup of a restriction semigroup $S$ is {\em full\/} if it contains all of the projections of $S$. 
A homomorphism $S \longrightarrow T$ of restriction semigroups is $P$-{\em separating\/} (or projection-separating) if it is injective on $P_S$.   

Likewise, the term `congruence' is used appropriate to context. We denote the greatest $P$-separating congruence by $\mu$ and observe that $\mu \subseteq \,\HP$.  
In the standard terminology, a congruence $\rho$ on a semigroup is {\em perfect\/} if $(a \rho)(b \rho) = (ab)\rho$. 

The {\em least monoid congruence\/} on a restriction semigroup $S$, denoted $\sigma$, is the least congruence that identifies all projections.  It is well known that
\[\sigma = \{(a,b) \in S: e a = e b \text{ for some } e \in P_S\} = \{(a,b) \in S: a f = bf \text{ for some } f \in P_S\}. \] It is clear that any principal order ideal of $S$ is contained in a $\sigma$-class and that any $\sigma$-class is an order ideal of $S$.

A restriction semigroup $S$ is {\em proper\/} if $\mathrel{\RP} \cap \mathrel{\sigma}\,  = \, \mathrel{ \LP} \cap \mathrel{\sigma} = \iota$ (where $\iota$ is the identical relation).   From this definition it is immediate that  $\sigma \cap \mu = \iota$, that is, $S$ is a subdirect product of $S/\sigma$ and $S/\mu$.  

Given a restriction semigroup $U$, a {\em cover for $U$ [over a monoid $T$]\/} is a proper restriction semigroup $S$ [such that $S/\sigma \cong T$] having $U$ as a $P$-separating homomorphic image.

We recall the primary definitions of this paper.  A restriction semigroup $S$ is {\em almost perfect\/} if it is proper and $\sigma$ is perfect.  A restriction monoid is {\em perfect\/} if it is proper, $\sigma$ is perfect and each $\sigma$-class has a greatest element.  By extension of the usual term for inverse semigroups, if a proper restriction semigroup has just the latter property, we term it {\em $F$-restriction\/}.  Since the projections form a $\sigma$-class, such a semigroup is of necessity a monoid.  In any such monoid, let $m_a$ denote the greatest element of $a\sigma$, so that $a \sigma = m_a \darrow$. 

The following result plays a key role in this paper (and motivates consideration of subhomomorphisms although, in the body of the paper, full consideration to the latter will only be given in Section~\ref{proper}).
A mapping $\alpha: S \longrightarrow T$ of restriction semigroups is a {\em subhomomorphism\/} if $(a \alpha) (b \alpha) \leq (ab) \alpha$ for all $a,b \in S$ (and the unary operations are respected, which if $S$ is a monoid, interpreted as either a plain monoid or a restriction monoid, and $\alpha$ is monoidal, occurs automatically). The reader should beware that, for inverse semigroups, Lawson \cite{lawson_book} calls such maps dual prehomomorphisms, while for Petrich \cite{petrich} a prehomomorphism is a subhomomorphism that, in addition, respects inversion (also see Section~\ref{inverse}).

 \begin{lemma}\label{sub} Let $M$ be a proper $F$-restriction monoid and put $T = M/\sigma$.
 Then the map 
 $\kappa: T \longrightarrow S$, $(a\sigma) \kappa = m_a$, is a monoidal subhomomorphism, which is a homomorphism if and only if $M$ is perfect.
 \end{lemma}
 {\bf Proof}. In these terms, the inclusion $(a \sigma) (b\sigma) \subseteq (ab) \sigma$ translates into the inequality $m_a m_b \leq m_{ab}$, with equality in one being equivalent to equality in the other.\eop\\

The motivation for the term $T$-proper comes from consideration of generators.  Suppose a restriction semigroup $S$ is generated, as such, by a subset $X$. It is well known that every non-projection of $S$ is expressible as the product of a projection and an element of the subsemigroup $\langle X \rangle $ generated by $A$.  Thus for every $s \in S$, either $s \in P_S$ or $s \leq t$ for some $t \in \langle X \rangle$, in which case, $s = s^+t = t s^*$. Let $T$ denote the `plain' submonoid of $S^1$ generated by $X$. Then $S = P_S T$.  In this case we say that $T$ {\em $P$-generates\/} $S^1$.

Finally, recall that any inverse semigroup $(S, \cdot, ^{-1})$ may be regarded as a restriction semigroup by setting $x^+ = xx^{-1}$ and $ x^*= x^{-1}x$ and `forgetting' the inverse operation. In that case, $P_S = E_S$. 
Now $\sigma$ is the least group congruence and the term `$E$-unitary' is more commonly used, rather than `proper'.

The inverse semigroups of most significance in this paper are the topic of the next subsection.
\subsection{Munn semigroups and the (generalized) Munn representations}\label{M}
If $Y$ is a semilattice,  $T_Y$  denotes the {\em Munn\/} semigroup on $Y$: the inverse subsemigroup of the symmetric inverse semigroup ${\cal I}_Y$ consisting of the isomorphisms between {\em principal\/} ideals of $Y$. For an exposition of this semigroup, its basic properties, and application to inverse semigroup theory, see \cite{howie}. The dual semigroup $T_Y^r$ again consists of the isomorphisms between principal ideals, but with functions instead written on the left and composition reversed.  Denote by $\gamma$ the anti-isomorphism $T_Y \longrightarrow T_Y^r$ that is induced by inversion.

Less familiar is the inverse subsemigroup of ${\cal I}_Y$, which we denote $TI_Y$, consisting of the isomorphisms between {\em arbitrary\/} ideals of $Y$, which plays a central role in this paper. (In $\cite{petrich}$, this semigroup is denoted $\Sigma(Y)$.)  In Proposition~\ref{C(T_Y)} we shall provide an alternative representation of  $TI_Y$ (and cite another one at the end of the following subsection).

Throughout this paper, the domain and range of a partial mapping $\alpha$ are denoted $\Delta \alpha$ and $\nabla \alpha$.

Let $S$ be a restriction semigroup and put $Y = P_S$.  The {\em generalized (right) Munn representation\/} $\theta: S \longrightarrow T_Y$ is defined by $a \mapsto \theta_a$, where $\Delta \theta_a = a^+ \darrow$ and, for $ e \leq a^+$, $e \theta_a = (ea)^*$.  In the case of an inverse semigroup, this reduces to the usual Munn representation.  The generalized left Munn representation is the dual map $\psi: S \longrightarrow T_Y^r$, defined by $a \mapsto \psi_a$, where $\Delta \psi_a = a^* \darrow$ and, for $ f \leq a^*$, $f \psi_a = (af)^+$.  For each $a$,  $\theta_a$ and $\psi_a$ are mutually inverse isomorphisms between their respective domains. In the sequel, we shall omit the qualifier `generalized'.

The bulk of the following result is in \cite[Proposition 5.2]{FGG1}.  A broad generalization, framed in the language used herein, was found by the author \cite{jones_P-restriction}.

\begin{result}\label{representations} Let $S$ be a restriction semigroup and put $Y = P_S$.  The maps $\theta$ and $\psi$ are biunary, projection-separating homomorphisms from $S$ onto full subsemigroups of $T_Y$ and $T_Y^r$, respectively, related by $\psi = \theta \gamma$.  Each induces the greatest $P$-separating congruence $\mu$ on $S$.

If $S$ is a monoid, then $\theta$ and $\psi$ are monoidal.
\end{result}


\subsection{The monoids $C(S)$ and their representations }\label{C}
We refer the reader to \cite{S1} for more details of the basic results cited here. There Szendrei defined $C(S)$ for the case of restriction semigroups in general, closely following \cite {GS} (itself based on the one-sided notion of El Qallali \cite{EQ} and extending the definition in the case of inverse semigroups \cite[V.2.2]{petrich}).

A {\em permissible\/} subset $A$ of a restriction semigroup $S$ is a nonempty order ideal such that $a^+b = b^+a$ and $ab^* = ba^*$ for all $a,b \in A$.  Note that if $S$ is proper, then the latter criteria are each equivalent to the property that all members of $A$ are $\sigma$-related.  We summarize the properties that will be used in the sequel. Originally proved in \cite{GS} for the subclass of `weakly ample' semigroups, the proofs carry over immediately to all restriction semigroups.

\begin{result}\cite[Theorems 3.1, 3.2]{S1}, \cite[Proposition 3.8]{GS} \label{C(S)} Let $S$ be a restriction semigroup.  The set $C(S)$ of all permissible subsets of $S$ is a restriction monoid, under multiplication of subsets, with identity $P_S$, where if $A \in C(S)$, then $A^+ = \{a^+ : a \in A\}$ and $A^* = \{a^* : a \in A\}$; its natural partial order is inclusion; its semilattice of projections consists of the ideals of $P_S$, under inclusion (and is thus $C(P_S)$).  The map $\tau_S: a \mapsto a \darrow$ embeds $S$ in $C(S)$; if $S$ is a monoid, then the embedding is monoidal.

The monoid $C(S)$ is proper if and only if $S$ is proper.  In that event, if $A \in C(S)$ then $A \subseteq a\sigma_S$ for any $a \in A$ and $a\sigma_S$ is the greatest element of $A \sigma_{C(S)}$, so that $C(S)$ is an $F$-restriction monoid. The monoids $S/\sigma$ and $C(S)/\sigma$ are therefore isomorphic.

Let $\theta: S \longrightarrow T$ be a homomorphism of restriction semigroups whose image is an order ideal.  Then the mapping $\hat{\theta}: C(S) \longrightarrow C(T)$, defined by $A \mapsto A \theta$, is a monoidal homomorphism such that $\tau_S \, \hat{\theta} = \theta \, \tau_T$.
\end{result}

In the event that $S$ is proper and we put $T = S/\sigma$, then the map $\kappa: T \longrightarrow C(S)$ defined in Lemma~\ref{sub} is just the injection map.

\begin{propo}\label{C(T_Y)}  Let $Y$ be a semilattice. The map $\Sigma:  C(T_Y) \longrightarrow TI_Y$, where $A\Sigma$ is the union of the members of $A$ (regarded as relations on $Y$), is an isomorphism such that $\tau_{T_Y} \Sigma$ is the inclusion map $T_Y \longrightarrow TI_Y$.
\end{propo}

{\bf Proof}. If $A\in C(T_Y)$, then $A$ is an order ideal that further satisfies, in the context of inverse semigroups, $\alpha\alpha^{-1} \beta = \beta \beta^{-1} \alpha$ for all $\alpha, \beta \in A$.  The latter property says that any two members of $A$ agree on the intersection of their domains, and so $A \Sigma$ is a well-defined order isomorphism between the ideals consisting of the unions of the domains and ranges, respectively, of the members of $A$.

Conversely, if $\alpha \in  TI_Y$, let $A_\alpha = \{ \alpha_{\vert\, e\darrow} : e \in \Delta \alpha\}$.  It is now routinely verified that $A_ \alpha \in C(T_Y)$ and that the maps $A \mapsto A_\alpha$ and $\Sigma$ are mutually inverse isomorphisms.

If $\alpha \in T_Y$, then  since $\alpha \darrow$ is a principal order ideal, the union of its members is just $\alpha$ itself, yielding the final statement.\eop

\begin{propo}\label{Munn extension} For any restriction semigroup $S$, the Munn representation $\theta: S \longrightarrow T_{P_S}$ induces a homomorphism $\otheta: C(S) \longrightarrow TI_{P_S}$ such that $\tau_S \otheta = \theta$.  For each $A \in C(S)$, $\otheta_A : A^+ \longrightarrow A^*$ and, if $e \in A^+$, $e \otheta_A = b^*$, where $b \in A$, $b^+ = e$.

\end{propo}

{\bf Proof}. For convenience, put $Y= P_S$.  Since the image of $\theta$ in $T_Y$ is full (Result~\ref{representations}), it is an order ideal.  By Result~\ref{C(S)}, $\hat{\theta}: C(S) \longrightarrow C(T_Y)$ is a well defined homomorphism. Then $\otheta = \hat{\theta} \Sigma$ is a homomorphism. If $A \in C(S)$, then  $A \otheta = \bigcup_{a \in A} \theta_a$, from which the stated formulas follow (putting $b = ea$, for any $a\in A$ such that $e \leq a^+$).  Combining the two previous results, $\tau_S \otheta = \tau_S\, \hat{\theta} \Sigma = \theta \,\tau_{T_Y}\, \Sigma = \theta$.\eop\\

We will call $\otheta$ the {\em extension\/} of $\theta$ to $C(S)$.
It is easily checked that if $S$ is a monoid, then all the homomorphisms above preserve the relevant identity elements.
The combination of this proposition with Lemma~\ref{sub} yields the following corollary, which is at the heart of this paper.

\begin{corol}\label{T-rep} Let $S$ be a proper restriction semigroup. Put $T = S/\sigma$ and $Y = P_S$. Then the composition $\kappa \otheta: T \longrightarrow TI_Y$ is a monoidal subhomomorphism, which is a homomorphism if and only if $\sigma$ is perfect on $S$ (that is, $S$ is almost perfect).
\end{corol}

\begin{figure}[h]
\begin{center}
$\begin{CD}
T  @ > \kappa >>  C(S) @ > \hat{\theta} >> C(T_Y)                             @>\Sigma >> TI_Y\\
    @.                                 @A\tau_S A A                @A \tau_{T_Y} AA  @.                              \\
    @.                       S       @> \theta >>             T_Y                                   @.
\end{CD}$
\end{center}
\caption{The mappings in this section.}\label{figure}
\end{figure}

We state without proof, since it is not used in this paper, that for any semilattice $Y$, $TI_Y$ is also isomorphic to the Munn semigroup of the semilattice $C(Y)$; and the homomorphisms $\hat{\theta}: C(S) \longrightarrow C(T_{P_S})$ and $\otheta: C(S) \longrightarrow TI_Y$ are equivalent to the Munn representations of $C(S)$.

\section{$T$-properness and almost $T$-properness}\label{properness}

Recall from Section~\ref{intro} that a restriction monoid $M$, with submonoid $T$, is $T$-{\em proper\/} if $M = P_M T$ and $\mathrel{\sigma}$ separates $T$; and that a restriction semigroup $S$ is {\em almost $T$-proper\/} if $C(S)$ is $T$-proper (with respect to a submonoid $T$).

The equation $M = P_M T$ merely states that $T$ $P$-generates $M$, as a restriction monoid (in the terminology introduced in Section~\ref{prelim}). By abuse of terminology, given any monoid $T$, we may call $M$ $T$-proper when it is $T'$-proper with respect to a submonoid $T'$ isomorphic to $T$.

In Propositions~\ref{T-proper} and \ref{almost T-proper}, we show that these definitions are equivalent, respectively, to those of perfection and almost perfection. Each form has its own benefits: the `$T$-proper' versions have the virtue of being consistent with prior analogous use of the term for Ehresmann semigroups \cite{GG2} and allow ready identification of examples; the `perfect' versions have the virtue of independence from $T$ and of straightforward verification once $\sigma$ has been computed.

We begin with element-wise criteria for the monoidal case, illustrating calculations prevalent in the sequel. The global characterizations immediately follow.

\begin{lemma}\label{basic} A restriction monoid $M$ is $T$-proper if and only if for all $m \in M$, $m \leq t$ for some unique $t \in T$; equivalently, $m = m^+ t$ for some unique $t \in T$; and equivalently if and only if for all $m \in M$, $m = t m^*$ for a unique $t \in T$.  In that event, $M$ is necessarily proper and $T \cong M/\sigma$.
\end{lemma}

{\bf Proof}. It was noted in Section~\ref{prelim} that $M$ is $P$-generated by $T$ if and only if for all $m \in M$, $m \leq t$ for some $t \in T$.  Let $t, u \in T$.  If $t \mathrel{\sigma} u$, then $et = eu = a$, say, for some $e \in P_S$, so that $a \leq t, u$ and $t = u$.  Conversely, if $a = et = eu$ then $t \mathrel{\sigma} u$. The equivalence of the next two statements follows from the discussion of $\sigma$ in Section~\ref{prelim}.

To prove the last statements, assume $M$ is $T$-proper and suppose that $(a, b) \in \, \mathrel{\sigma} \cap \RP$.  Write $a = a^+t$ and $b = b^+u$, for $t, u \in T$. Then $t \mathrel{\sigma} a \mathrel{\sigma} b \mathrel{\sigma} u$, so that $t = u$; and $a^+ = b^+$, so $a = b$.  The dual statement is proved analogously.  Now the submonoid $T$ is a transversal of the $\sigma$ classes, whereby it is isomorphic to the quotient monoid.\eop
\begin{propo}\label{T-proper} The following are equivalent for a proper restriction monoid $M$:

\be[(i)]
\item $M$ is $T$-proper for some monoid $T$;

\item $M$ is $M/\sigma$-proper;

\item $M$ is perfect.

\ee
\end{propo}

{\bf Proof.}  The equivalence of (i) and (ii) follows from the last statement of the lemma.  If either holds, then each $\sigma$-class has a greatest element, again by the lemma, and these greatest elements are precisely the members of the submonoid $T$.  So, in the notation of Lemma~\ref{sub}, $m_a m_b = m_{ab}$ for all $a,b \in M$.  As noted in the proof of that lemma, this property is equivalent to perfection of $\sigma$.  Conversely, if $M$ is perfect, then it is $T$-proper with respect to $T = \{m_a: a \in M\}$.\eop\\

Turning now to an arbitrary restriction semigroup, we may of course consider $T$-properness of the monoid $M = S^1$. However, that is not general enough for our purposes, as will be seen in the sequel.

By Lemma~\ref{basic}, that a semigroup be almost $T$-proper is equivalent to the property that every $A \in C(S)$ be contained in a unique member of $T$. Also by that lemma, in that case $C(S)$ and, therefore, $S$ itself is necessarily proper. The next result is the analogue of Proposition~\ref{T-proper}. 

\begin{propo}\label{almost T-proper} The following are equivalent for a proper restriction semigroup $S$:
\be[(i)]
\item $S$ is almost $T$-proper for some monoid $T$;

\item $S$ is almost $S/\sigma$-proper;

\item $S$ is almost perfect.

\ee
\end{propo}

{\bf Proof}. According to Result~\ref{C(S)}, $C(S)/\sigma \cong S/\sigma$.  The equivalence of (i) and (ii) then follows from the definition and from applying the corresponding equivalence in Proposition~\ref{T-proper} to $C(S)$.

Again by Result~\ref{C(S)}, the $\mathrel{\sigma}$-classes of $S$ are already the greatest elements of the $\mathrel{\sigma}$-classes of the monoid $C(S)$. Applying Proposition~\ref{T-proper}, $C(S)$ is therefore $C(S)/\sigma$-proper if and only if these greatest elements form a submonoid under multiplication of subsets.  That is just the statement that $\mathrel{\sigma}$ itself is a perfect congruence.\eop\\

In the case of a monoid, the relationship between these two definitions needs to be elucidated.  We return to the original terminology of `perfection'.

\begin{propo}\label{T and almost T}  (1)  Every perfect restriction monoid is almost perfect; the converse need not hold.  (2)  If $S$ is a proper restriction semigroup without identity and $S^1$ is almost perfect, then $S$ is almost perfect; the converse holds if 
$S/\sigma$ has trivial group of units.

\end{propo}

{\bf Proof}. The direct statement in (1) is obvious from the definitions.  Example~\ref{example} shows that the converse need not be true.

To prove the direct statement in (2), note that in this case $\mathrel{\sigma}$ on $S$ is simply the restriction of that on $S^1$. Now by assumption $\mathrel{\sigma}$ is perfect on $S^1$. But the only product in $S^1$ that yields 1 is $1 \cdot 1$, so $\mathrel{\sigma}$ is also perfect on $S$ and $S$ is therefore also almost perfect.  For the converse, if $\mathrel{\sigma}$ on $S$ is perfect, then so is $\mathrel{\sigma}$ on $S^1$, under the stated assumption.\eop\\

Two natural motivating classes of examples for $T$-properness and almost $T$-properness are (i)  free (and certain relatively free) restriction semigroups and monoids, and (ii) factorizable and almost factorizable restriction semigroups (in one - and two - sided versions).  Rather than break the flow at this point, we refer the reader to Sections~\ref{examples} and \ref{embedding theorem}, respectively, for the precise connection.
The constructions in Sections~\ref{covers} and \ref{X} provide a ready source of further examples.

\section{Almost perfect covers}\label{covers}
We present in Theorem~\ref{newcovers} a far-reaching generalization $S_{T, R}$ of a construction first explicitly stated in \cite[Theorem 9.1]{jones_P-varieties}, although special cases have appeared elsewhere in the literature.  Its simplicity is somewhat obscured by the generality in which we frame it.  The pairs $(S, R)$ of restriction semigroups of most interest satisfy $S  \leq R \leq C(S)$ (identifying $S$ with its image in $C(S)$ under $\tau_S$). In particular, the pairs $(S, C(S))$ and $(S, S^1)$ play distinct major roles in this paper, as demonstrated in Corollaries~\ref{C(S) cover} and \ref{T-proper cover} respectively.  
The slightly more abstract setting that we have chosen simplifies the notation considerably and specializes more straightforwardly to the monoidal setting.

Because the specific monoid $T$ is crucial to this section, we tend to prefer the `$T$-proper' terminology instead.

Let $S$ be a restriction semigroup and $R$ a restriction monoid that contains $S$, such that for each $r \in R$ the following conditions and their duals are satisfied: (i) there exists $e \in P_S$, $e \leq r^+$, and (ii) for any such $e$, $er \in S$. Clearly this criterion is satisfied if $R = S^1$.  It will be shown after the theorem that this is also the case for the other instances described in the previous paragraph.  Let $T$ be a `plain' monoid and $\alpha: T \longrightarrow R$ a monoidal homomorphism. (In Section~\ref{proper}, $\alpha$ will be allowed to be a subhomomorphism).  Let 
\[S_{T,R} = \{(a,t ) \in S \times T : a \leq t \alpha \text{  in  } R\}.\]
Write $S_T$ in case $R = S^1$. 
 
\begin{theorem} \label{newcovers} As just defined, $S_{T,R}$ is an almost $T$-proper restriction subsemigroup of $S \times T$, with $S_{T,R} /\sigma \cong T$, now regarding $T$ as a restriction monoid,  and the first projection is a $P$-separating homomorphism onto  a full subsemigroup of $S$. Further:

\be

\item if $S \subseteq (T \alpha)\!\darrow $,  then $S_{T, R}$ is a subdirect product of $S$ and $T $ and is therefore a cover of $S$ over $T$;  in particular, if $R = S^1$ and $T\alpha$ $P$-generates $R$, $S_T$ is a cover of $S$;

\item if $S$ is a monoid and $\alpha$ is monoidal, then $S_T$ is $T$-proper.
\ee
\end{theorem}

{\bf Proof}.  Let $(a,t), (b,u) \in S_{T,R}$.  By compatibility of the natural partial order on $R$ and the homomorphism property of $\alpha$, $ab \leq (t \alpha)(u\alpha) \leq (tu) \alpha$, so 
$S_{T,R}$ is closed.  (Observe that $\alpha$ only need be a subhomomorphism for this to hold.) For any $(a,t) \in S_{T,R}$, $(a, t)^+ = (a ^+, 1)$ and $(a,t)^* = (a^*, 1)$, where $a^+, a^* \leq 1 = 1 \alpha$, so $S_{T,R}$ is a restriction subsemigroup of $S \times T$.  The semilattice $P_{S_{T,R}} = \{(e,1) : e \in P_S \} \cong P_S$.

By the assumption on $S$ and $R$, for any $t \in T$ there exists $e \in P_S$ such that $e \leq t \alpha$ and $e (t\alpha) \in S$, so that $(e (t \alpha),t) \in S_{T,R}$.  That is, the second projection map $S_{T,R} \longrightarrow T$ is surjective.  In general, this is not true of the first projection map $S_{T,R} \longrightarrow S$. However since $(e,1) \mapsto e$, it is projection-separating and its image is a full subsemigroup.

Next we compute $\sigma$.  Let $(a,t), (b,u) \in S_{T,R}$. If $(a,t) \mathrel{\sigma} (b,u)$ then from the triviality of $\sigma$ on $T$, $t = u$. Conversely, if $t = u$, then $a, b \leq t\alpha$.  In the case $t = 1$, then $(a,1), (b,1)$ are both projections and so are $\sigma$-related; otherwise, $(b^+, 1)(a,t)  =
 (b^+ a^+(t\alpha), t) = (a^+b^+ (t \alpha), t) = (a^+, 1)( b, t)$, so that $(a,t) \mathrel{\sigma} (b,t)$.  Therefore the $\mathrel{\sigma}$-classes of $S_{T,R}$ are the sets $S_t = \{(a,t): a \leq t \alpha\}$, $t \in T$, and $S_{T,R}/\sigma \cong T$.

Now suppose $((a,t) , (b,u)) \in \, \RP\! \cap\, \sigma$ in $S_{T,R}$.  Then $t = u$ and $a = a^+ t = b^+ u = b$. In combination with the dual argument, this shows that $S_{T,R}$ is proper.

To prove almost $T$-properness, using Proposition~\ref{almost T-proper} it must be shown 
 that $S_t S_u = S_{tu}$. Clearly $S_t S_u \subseteq S_{tu}$; conversely, suppose $(c,tu) \in S_{tu}$, where $c \leq (tu) \alpha = (t\alpha)(u\alpha)$ (using the homomorphism property of $\alpha$) and $t,u \in T$.  Then $c = c^+(t \alpha) ( u \alpha)c^*$, where $c^+, c^* \in P_S$, so that $c^+(t\alpha),  (u \alpha) c^*\in S$, using the assumption (ii) on $S$ and its dual. Therefore $(c, tu) = ( c^+(t \alpha), t)( ( u \alpha)c^*, u) \in S_t S_u$.
 

To prove 1, suppose that $S \subseteq (T\alpha)\darrow$.  Then for each $a \in S$, $(a,t) \in S_{T,R}$ for some $t \in T$.
Thus $S_{T,R}$ is a subdirect product of $S$ and $T$ and the first projection map, being already projection-separating, is a covering map.

To prove 2, suppose that $S$ is a monoid (and $R = S$). Since $S_T$ contains $(1,1) $, it is a monoid. Now for each $t \in T$, $(t \alpha, t )$ is the greatest element in the $\sigma$-class $S_t $. Thus $S_T$ is $F$-restriction and therefore perfect (and so $T$-proper).\eop\\

We detail the special case $R = S^1$ since for each restriction semigroup it produces a very simple almost perfect cover.  In Corollary~\ref{T-proper cover alt} we will provide a simple alternative representation of the coverings in this corollary by means of the $W$-semigroup construction. It will be shown below that, in the monoidal case, all perfect covers may be found as in the corollary.  That when applied to an {\em inverse\/} semigroup $S$, this covering does not in general produce another inverse semigroup is at the heart of the divergence of our work from inverse semigroup theory (see Section~\ref{inverse}).

\begin{corol} \label{T-proper cover} Let $S$ be a restriction semigroup. If $T$ is a monoid and $\alpha: T \longrightarrow S^1$ is a monoidal homomorphism, the image of which $P$-generates $S^1$, then $S_T = \{(s,t) \in S \times T: s \leq t\alpha \text{  in  } S^1\}$ is an almost $T$-proper cover of $S$ that is a subdirect product of $S$ and $T$.  In particular, if  $S$ is generated, as a restriction monoid, by a subset $X$ and $T$ is the `plain' submonoid of $S^1$ generated by $X$, then $S_T = \{(s,t) \in S \times T: s \leq t \text{  in  } S^1\}$ is such a cover.

If $S$ is a monoid to begin with, then the above covers are $T$-proper and monoidal.

\end{corol}

There is one immediate consequence that will be used later.

\begin{corol}\label{special cover} Any proper restriction semigroup [monoid] $S$ has an almost perfect [perfect] cover over $T = S/\sigma$ that belongs to the variety [of monoids] generated by $S$.

\end{corol}

{\bf Proof}. If $N$ is {\em any\/} proper cover of $S$ over $T$, then the congruence induced by the homomorphism upon $S$ is $P$-separating and so contained in $\mu$, and since $\mu \cap \sigma = \iota$ (see Section~\ref{prelim}), $N$ is a subdirect product of $S$ and $N/\sigma \cong S/ \sigma$.\eop\\

Another case of special interest occurs when $S$ is generated, as a restriction semigroup, by a subset $X$: the cover $S_{X^*}$ associated with the homomorphism $\alpha: X^* \longrightarrow S^1$.  In this case the cover is ample.\\

We now present another useful specialization. In Theorem~\ref{P-sep covers} and its corollary, it will be shown that it is in a sense as general as the original theorem. In the following, we shall identify a restriction semigroup $S$ with its image $S\tau_S$ in $C(S)$, according to Result~\ref{C(S)}.  However we shall at times revert to the use of $\tau_S$ when clarification is needed.  Recall that the partial order in $C(S)$ is inclusion, whereby the definition of $S_{T,R}$ below is the specialization of that in Theorem~\ref{newcovers}.  Note that Corollary~\ref{T-proper cover} represents the simplest case of this result, since by Result~\ref{C(S)}, the embedding of $S^1$ in $C(S)$ is monoidal, so that the element 1 can be identified with the identity of $C(S)$.

\begin{corol}\label{C(S) cover} Let $S$ be a restriction semigroup, $R$ a submonoid of $C(S)$ that contains $S$, $T$ a `plain' monoid and $\alpha: T \longrightarrow R$ a monoidal homomorphism. The almost $T$-proper semigroup
\[S_{T,R} = \{(a,t ) \in S \times T : a \in t \alpha \}\]
is well defined.  In particular, $S_{T, C(S)}$ is well defined.

For any semilattice $Y$, `plain' monoid $T$ and monoidal homomorphism $\alpha: T \longrightarrow TI_Y$, the semigroup $(T_Y)_{T, TI_Y}$ is a well defined, almost $T$-proper restriction semigroup.  If $Y$ has an identity element and $\alpha: T \longrightarrow T_Y$, then $(T_Y)_T$ is a $T$-proper monoid.
\end{corol}

{\bf Proof}. We verify that the criteria (i) and (ii) are met.  For the sake of clarity, we set aside the identification of $S$ with $S \tau_S$.  
Let $A \in R$.  Then $A$ is an order ideal of $S$ that satisfies $a^+b = b^+a$ for all $a,b \in A$. Here $A^+$ is an ideal of $P_S$ and so contains $e \darrow = e \tau_S$, for some $e \in P_S$.  For any such $e$, there exists $b \in A$ such that $b^+ = e$. Then $e\darrow \! A =  b\darrow = b\tau_S$: for if $f \leq e$ and $a \in A$, then $fa = f b^+a = f a^+b \in b \darrow$; and if $c \in b \darrow$, then $c \in A$ and $c = b^+c  \in e\darrow \! A$. The dual statement follows from the self-duality of the pair $(S, R)$.

The statements in the second paragraph follow from the identification of $C(T_Y)$ with $TI_Y$ in 
 Proposition~\ref{C(T_Y)}. (Once (i) and (ii) have been established, Theorem~\ref{newcovers} could itself be applied to the pair ($T_Y, TI_Y$).)\eop\\ 

The following theorem may be regarded as a converse of Theorem~\ref{newcovers}. Its generality allows two distinct important applications. Again we identify $S$ with its image in $C(S)$ under $\tau_S$, except where additional clarity is required. Refer to Result~\ref{C(S)} for the relevant properties of $C(S)$.

\begin{theorem}\label{P-sep covers} Let $N$ and $S$ be restriction semigroups, with $T = N/\sigma$.  If $N$ is almost perfect and $\beta: N \longrightarrow S$ is a $P$-separating homomorphism whose image is full in $S$, then $N \cong S_{T,C(S)}$, with respect to $\alpha = \kappa \hat{\beta}: T \longrightarrow C(S)$.

Let $N$ and $M$ be restriction monoids. If $N$ is perfect and $\beta: N \longrightarrow M$ is a $P$-separating homomorphism whose image is full in $M$, then $N \cong M_T$, where
$\alpha = \kappa \beta: T \longrightarrow M$.
\end{theorem}

{\bf Proof}.  By Theorem~\ref{C(S) cover}, $S_{T,C(S)}$ is well defined.  Recall that $\kappa: T \longrightarrow C(N)$ and $\hat{\beta}: C(N) \longrightarrow C(S)$ is the monoidal extension of $\beta$, according to Result~\ref{C(S)}, noting that since $N \beta$ is full in $S$, it is an order ideal.  See Figure~\ref{figure 2}.  Since $1 \kappa = P_N = 1_{C(N)}$, $\alpha$ is monoidal.
\begin{figure}[h]
\begin{center}
$\begin{CD}
T  @ > \kappa >>  C(N) @ > \hat{\beta} >> C(S)\\
    @.                                 @A\tau_N A A                @A \tau_S AA\\
    @.                       N       @> \beta >>             S
\end{CD}$
\end{center}
\caption{The mappings in the proof of Theorem~\ref{P-sep covers}.}\label{figure 2}
\end{figure}


Define $\omega: N \longrightarrow S_{T,C(S)}$ by $n \omega = (n\beta, n \sigma)$.  Clearly $\omega$ is a homomorphism.  Since $\beta$ is $P$-separating, so is the congruence on $N$ that it induces.  Since $N$ is proper, $\mu \cap \sigma = \iota$, so $\omega$ is injective.

Let $(s,t) \in S_{T,C(S)}$, that is, $s  \leq t\alpha$ or, more precisely, $s \darrow\, \subseteq t \alpha$ in $C(S)$, so that $s \in t \alpha$. Therefore $s  = a \beta$ for some $a $ in the $\sigma$-class $t = a\sigma$ of $N$.  But then $(s,t) = a \omega$. So $\omega$ is surjective.


In the monoidal case, $\kappa $ maps directly into $N$ and $\omega$ may now be regarded as mapping directly into $M_T$; the proof then proceeds almost identically.\eop\\

The first application of this theorem is a description of the almost perfect covers of a restriction semigroup.

\begin{corol}\label{all covers} Let $N$ be an almost perfect cover of a restriction semigroup $S$, via the homomorphism $\beta$. Put $T = N/\sigma$. Then $N \cong S_{T, C(S)}$, where $\alpha = \kappa \hat{\beta}: T \longrightarrow C(S)$ and $S \subseteq (T \alpha)\!\darrow$ is satisfied (cf 1 in Theorem~\ref{newcovers}).

If $N$ is a perfect, monoidal cover of the restriction monoid $M$, again via $\beta$, then $N \cong M_T$, where
$\alpha = \kappa \beta: T \longrightarrow M$ and $T\alpha$ $P$-generates $M$. 
\end{corol}

{\bf Proof}. Here $\beta$ is, further to the theorem above, surjective.  Let $s \in S$, $s = a \beta$, say.  Reversing the argument in the proof of surjectivity above, $s \leq t \alpha$, where $t = a \sigma \in T$.  The monoidal case proceeds similarly.\eop\\

In the case that $S$ itself is almost perfect, regarded as its own cover, that is, $\beta$ is the identity map, put $T = S/\sigma$.  Then $S_{T,C(S)} \cong S$, since if $(s,t) \in S_{T,C(S)}$, then $t$ must be $s\sigma$.  Likewise, if $M$ is a perfect monoid, then $N \cong M_T$. 

Alternatively, $S$ may also be represented in the form $U_{T,R}$ via its Munn representation, as follows.  This result is a sort of precursor to the $W$-semigroup representation Theorem~\ref{converse}.  See the remarks following that theorem.

\begin{corol}\label{R-rep}  Let $S$ be an almost perfect restriction semigroup, with $T = S/\sigma$ and $Y = P_S$.  Then $S \cong F_{T,C(F)}$, where $F \cong S/\mu$ is the image of $S$ in $T_Y$ under the Munn representation.
\end{corol}


\section{The general $W$-product}\label{X}

Let $T$ be a monoid, $Y$ a semilattice and suppose there is a monoidal homomorphism $\alpha : T \longrightarrow TI_Y$, where (see Section~\ref{prelim}) $TI_Y$ is the inverse semigroup of isomorphisms between ideals of $Y$.  Adapting the usual language of actions, we say that $T$ acts on $Y$ (on the right) by {\em isomorphisms between ideals\/}.  If the image lies in $T_Y$, then we say that the action is by isomorphisms between {\em principal\/} ideals. Once more, in Section~\ref{proper} we will allow {\em sub\/}homomorphisms in the construction.  The relationship between our construction and the original $W$-product will be elucidated following Theorem~\ref{X semigroup}.

For $t \in T$, write $\alpha_t$ instead of $t\alpha$ and denote by $\Delta t$ and $\nabla t$, respectively, its domain and range.  In the usual notation of (partial) right actions, then, for $e \in \Delta t$, write $e^t$ for $e \alpha _t$.  

Consider the set
\[W(T, Y) = \{(t,f) \in T \times Y : f \in \nabla t\}  =  \{(t,e^t) \in T \times Y : e \in \Delta t\}.\]

The alternative form $(t,e^t)$ for $(t,f)$ results from the bijectivity of $\alpha_t$.  Each form will prove to be the more convenient one at different points.  The latter description reduces precisely to the original definition of $W(T,Y)$ (see Sections~\ref{intro} and \ref{X semigroup}) in the case that the action of each $t$ is fully defined.

The product is defined by:
\[ (t,e^t)(u, f^u) = (tu, (e^t f)^u).\]
Since $\Delta u$ is an ideal containing $f$ and $e^t f \leq f $, $(e^t f)^u$ is defined; likewise, $e^t f \in \nabla t$, $e^t f = g^t$, say.  Since the action is induced by a homomorphism, $g \in \Delta tu$ and $(e^t f)^u = (g^t)^u = g^{tu} \in \nabla tu$. Thus the operation $W(T,Y)$ is well defined. (Note that $\alpha$ only need be a subhomomorphism for this to hold.)

Unary operations are defined by:
\[(t,e^t)^+ = (1,e) \quad \text{and} \quad (t,f)^* = (1,f).\]

\begin{theorem}\label{X semigroup} Let $T$ be a monoid, $Y$ a semilattice, and $\alpha : T \longrightarrow TI_Y$ a monoidal homomorphism, that is, $T$ acts on $Y$ by isomorphisms between ideals.  Then $W = W(T,Y)$ is isomorphic to the semigroup $(T_Y)_{T, TI_Y}$ and is therefore 
an almost perfect restriction semigroup, with $P _W \cong Y$ and $W /\sigma \cong T $.

Further if $Y$ is a monoid and $\alpha$ is a homomorphism into $T_Y$, that is, $T$ acts on $Y$ by isomorphisms between principal ideals, then $W$ is isomorphic to $(T_Y)_T$ and is therefore a perfect monoid.
\end{theorem}

{\bf Proof}. The general conclusions will follow once the isomorphism is proved.  
Corollary~\ref{C(S) cover} shows that $(T_Y)_{T, TI_Y}$ is well defined.
 Recall that $Y \cong P_{T_Y} = E_{T_Y}$, under the map $e \mapsto e\darrow$.

Define $\Pi : (T_Y)_{T, TI_Y} \longrightarrow W = W(T,Y)$ by \[(\beta, t) \Pi = (t, e^t), \text{ where } e \darrow = \Delta \beta   .\]

By assumption, $\beta \subseteq t \alpha$, so $e \in \Delta t$.
  Therefore $(\beta, t) \Pi $ is well defined.  On the other hand, for any $(t, e^t) \in W$, let $\beta$ be the restriction of $\alpha_t$ to $e \darrow$. Then $(t, e^t) = (\beta, t) \Pi$, so $\Pi$ is surjective.  
  
To prove that $\Pi$ respects the binary operation, let $(\beta, t), (\gamma, u) \in (T_Y)_{T, TI_Y}$, where $\Delta \beta = e \darrow$ and $\Delta \gamma = f \darrow$. 
Now $ (\beta \gamma, tu)\Pi = (tu, g^{tu})$, where $ g\darrow = \Delta (\beta \gamma) = (\nabla \beta \cap \Delta \gamma) \beta^{-1} $.  Here since $\beta \subseteq t\alpha$,  $\nabla \beta = (e \beta)\darrow = e \alpha_t \darrow = e^t \darrow$ and so $\nabla \beta \cap \Delta \gamma = e^t \darrow \cap f \darrow = (e^t f) \darrow$, whereby $g \beta = g^t = e^t f$.  As in the definition of the operation on $W$, $ (e^t f)^u = (g^t)^u = g^{tu}$ and so $(tu, g^{tu}) = (t, e^t) (u, f^u) $.

Finally, $(\beta, t)^+ \Pi = (\beta^+, 1) \Pi = (1, e) = (t, e^t)^+$, since $\beta^+ $ is the identity map on $e\! \darrow$.  Likewise, $(\beta, t)^* \Pi = (\beta^*, 1) \Pi = (1, h)$, where $h\! \darrow = \nabla \beta = (e\beta)\!\darrow = (e \alpha_t)\! \darrow = e^t \! \darrow$, that is $(1,h) = (1, e^t) = (t, e^t)^*$. 
So $\Pi$ is a biunary isomorphism.

The monoidal case proceeds similarly.\eop\\

As noted in Section~\ref{intro}, the original $W$-semigroup construction \cite{FG, GS, S1} corresponds precisely to the special case whereby $\Delta t  = Y$ for all $t \in T$, that is, the action is by {\em endomorphisms\/} of $Y$. In that case, $W(T, Y) = T \times Y$ and so is a `reverse' semidirect product. (Observe that the action is not simply by endomorphisms, however, since these endomorphisms must be injective and their images must be ideals of $Y$).  The original construction of course includes the case that the action be by {\em automorphisms\/} of $Y$.

Our construction has a natural self-duality, due to the symmetry of $TI_Y$ under the anti-isomorphism $\delta \mapsto \delta^{-1}$.  Thus one can just as easily specialize to the case where $\nabla t = Y$ for all $t \in T$, that is, each map is surjective.  In the context of the original $W$-product, the dual construction must be used.

Proposition~\ref{factorizable} specializes Theorem~\ref{converse} below to the original construction.\\


The following are consequences of the definitions of the unary operations and the elementary properties of $S_{T,R}$ in Section~\ref{covers}.

\begin{lemma}\label{X Greens} Under the hypotheses of the theorem, if the pairs $ (t, e^t) = (t,g)$ and $ (u, f^u) = (u, h)$ belong to $ W = W(T,Y)$, then: $ (t, e^t) \RP (u, f^u)$ if and only if $e = f$;  $(t, g) \LP (u, h)$ if and only if $g = h$; $(t, g) \leq (u, h)$ if and only if $t = u$ and $g \leq h$;  $(t, g) \mathrel{\sigma} (u, h)$ if and only if $t = u$.

\end{lemma}

Using Lemma~\ref{X Greens}, the $\mathrel{\sigma}$-classes of $W$ are precisely the sets $A_t = \{(t, h): h \in \nabla t\} = \{(t, e^t): e \in \Delta t\}$, $t \in T$.  
Translating into this language the corresponding statement from (the proof of) Theorem~\ref{newcovers}, $A_t A_u = A_{tu}$ for all $t,u \in T$.  That is, in the original terminology, the $T$-properness of $C(W)$ is witnessed by the submonoid $T \kappa = \{A_t: t \in T\}$.

Likewise, in the case that $Y$ has an identity element and $T$ acts by isomorphisms between principal ideals, then $T$-properness of $W(T,Y)$ is witnessed by the submonoid $T\kappa = \{(t, g^t) : t \in T\}$ of $W$ itself, where for each $t$, $g$ generates the principal ideal $ \Delta t$ of $Y$.

\begin{propo}\label{X action} Under the hypotheses of the theorem, identify $P_W$ with $Y$ under $(e,1) \mapsto e$. Recall the homomorphism $ \kappa \otheta $ introduced in Corollary~\ref{T-rep}: the 
composition of the injection of $T$ in $C(W)$ with the extension of the Munn representation of $W$ to $C(W)$. This homomorphism is precisely $\alpha$ and so the original action of $T$ on $Y$ is equivalent to that induced by $\kappa \otheta$.

In the case where $Y = Y^1$ and the original action is by isomorphisms between principal ideals, the restriction of the Munn representation of $W$ itself to (the image under $\kappa$ of) $T$ is $\alpha$ and so induces an action of $T$ on $Y$ that is equivalent to the original one.
\end{propo}

{\bf Proof}. In the general case, see Figure~\ref{figure} for the mappings involved, with $W$ replacing $S$. If $t \in T$, then $t \kappa = A_t$, in the notation above.  According to Proposition~\ref{Munn extension}, $\Delta \otheta_{A_t} = A_t^+ = \{(1,e) : e \in \Delta t\}$ and, for $(1,e)$ in this domain,  $(1,e) \otheta_{A_t} = b^*$, where $b = (t,f^t) \in A_t$ and $(1,f) = b^+ = (1,e)$.  That is, $(1,e) (\kappa\otheta)_t = (t,e^t)^*   = (1, e^t) = (1, e \alpha_t)$.  Identifying $(1,e)$ with $e$, $\alpha$ and $\kappa \otheta$ therefore have the same domains and take the same values.
 
In the monoidal case,  if $t \in T$ and $\Delta t = g \darrow$, as above, then $\Delta \theta_ {(t, g^t)} = (t, g^t)^+ \darrow = (1,g)\darrow$.  For $e \leq g$, $(1,e) \theta_{(t,g^t)} = ((1,e)(t, g^t))^* = (t, e^t)^* = (1, e^t) = (1, e \alpha)$.  Again identifying $(1,e)$ with $e$, and $T$ with its image $T \kappa$ in $W$, the restriction to $T$ of the Munn representation agrees with $\alpha$.\eop

\section{Representing proper restriction semigroups as $W$-products}\label{represent}

Proposition~\ref{X action} suggests a straightforward route to the converse of Theorem~\ref{X semigroup}, via the homomorphism $\kappa \otheta$ introduced in Corollary~\ref{T-rep}, using Theorem~\ref{P-sep covers}.

\begin{theorem}\label{converse} Let $S$ be an almost perfect restriction semigroup.  Put $T = S/\sigma$ and $Y = P_S$.  Then $S \cong W(T,Y)$, where $Y = P_S$ and the action of $T$ on $Y$, by isomorphisms between ideals, is induced by the homomorphism $\kappa \otheta$, the composition of the injection of $T$ in $C(S)$ with the extension of the Munn representation of $S$ to $C(S)$.

If $S$ is a perfect restriction monoid, then $Y = Y^1$ and the action of $T$ on $Y$ is by isomorphisms between principal ideals,  induced by the Munn representation of $S$ itself.
\end{theorem}

{\bf Proof}. The Munn representation $\theta: S \longrightarrow T_Y$ is $P$-separating and its image is a full subsemigroup of $T_Y$, so by Theorem~\ref{P-sep covers}, $S \cong (T_Y)_{T,C(T_Y)}$, with respect to $\alpha = \kappa \hat{\theta}: T \longrightarrow C(T_Y)$.  Now by Proposition~\ref{C(T_Y)},  $(T_Y)_{T,C(T_Y)} \cong (T_Y)_{T, TI_Y}$, with respect to $\alpha \Sigma = \kappa \otheta$.  (Or the identifications in Corollary~\ref{C(S) cover} could have been used to combine those steps.)  Finally, by Theorem~\ref{X semigroup}, $(T_Y)_{T,TI_Y} \cong W(T, Y)$, with respect to the same map $\kappa \otheta$.

The proof in the monoid case proceeds similarly.\eop\\

Note that Corollary~\ref{R-rep} actually gave a sharper conclusion, based on the precise image $F$ of $S$ in $T_Y$, one that could be sharpened further by replacing $C(F)$ by the restriction monoid generated by $T\alpha$. The advantage of the $W(T,Y)$ formulation is that its parameters do not require specifying such subsemigroups.

In a different direction, observe that if $S$ is almost perfect, then $T$-properness of $C(S)$ implies that the latter semigroup is isomorphic to $W(T, P_{C(S)})$, via the action induced by its own Munn representation in $T_{P_{C(S)}} = T_{C(P_S)}$. At the end of Subsection~\ref{C}, it was remarked without proof that $T_{C(P_S)} \cong TI_{P_S}$
and that the Munn representation is equivalent to $\otheta$.

As promised before Corollary~\ref{T-proper cover}, we interpret the coverings presented there in terms of the $W$-product.  For any restriction semigroup $R$, with $Y = P_R$, the Munn representation $\theta: R \longrightarrow T_Y$ extends to a monoidal representation $R^1 \longrightarrow TI_Y$, by mapping $1$ to the identity of $TI_Y$.  Since this is in essence the Munn representation of $R^1$ in $T_{Y^1}$, we again denote it by $\theta$.  If $R$ is a monoid, this is just the original representation in $T_Y$.

\begin{corol} (Alternative version of Corollary~\ref{T-proper cover})\label{T-proper cover alt}  Let $S$ be a restriction semigroup, $T$ a monoid and $\alpha: T \longrightarrow S^1$ a monoidal homomorphism, the image of which $P$-generates $S$.  In the notation of Corollary~\ref{T-proper cover},  $S_T \cong W(T,Y)$, where $Y = P_S$ and the action between ideals of $Y$ is induced by $\alpha \theta$.

If $S$ is a monoid, then the action is by isomorphisms between principal ideals of $Y$, induced by the composition of $\alpha$ with the Munn representation of $S$.  In particular, if the plain submonoid $T$ of $S$  itself $P$-generates $S$, the action is induced by the restriction of the Munn representation to $T$.
\end{corol}

{\bf Proof}.  Here $S_T/\sigma \cong T$ and $P_{S_T} = P_S \times \{1\} \cong Y$. It only needs to be verified that the action exhibited in Theorem~\ref{converse} is equivalent to the one in the statement of this corollary.  For $t \in T$, the $\sigma$-class $t \kappa$ in $C(S_T)$ is $A_t = \{(a,t) : a \leq t \alpha \text{ in } S^1\}$.

We apply Proposition~\ref{Munn extension}. If $t \in T$ then \[\Delta t = \Delta_{A_t} = A_t^+ = \{(a,t)^+: a \leq t \alpha\} = \{(e,1): e \leq (t\alpha)^+\}.\] The last equation holds because if $a \leq t \alpha$, then $a^+ \leq (t 
\alpha)^+$; and if $e \leq (t \alpha)^+$, then $(e,1) = (e (t\alpha), t)^+$, where $e (t \alpha) \leq t \alpha$.

Now for $(e,1) \in \Delta t$, 
$(e,1)\otheta_{A_t} = (a,t)^*$, where $(a,t)^+ =  (e,1)$.  Since $a \leq t \alpha$, $a = e( t \alpha)$ and so 
$(e,1)\otheta_{A_t} = (e (t\alpha), t)^* =  ((e( t\alpha))^*, 1)$.

Identifying $P_S$ with $P_{S_T}$, the action of $t$ therefore has domain $P_{S_T}$, if $t = 1$, and domain $(t \alpha)^+ \darrow$ otherwise, with $e ^t = (e (t\alpha))^*$ in either case. The latter action is that induced by the composition of $\alpha$ with the (extension of the) Munn representation $\theta$.\eop\\

An explicit isomorphism $S_T \cong W(T,Y)$ is given by $(s,t) \mapsto (t, s^*)$, with inverse $(t, f) \mapsto ((t \alpha) f, t)$.

\section{Examples}\label{examples}
The $W$-product provides a simple mechanism for producing specific examples.  We begin this section with one promised in Proposition~\ref{T and almost T}.

 \begin{example}\label{example} An almost perfect restriction monoid need not be perfect.
\end{example}

{\bf Proof}. Let $Y = \ZZ$, under the reverse of the usual order. Let $T = \{x\}^*$, the free monoid on $\{x\}$, and let $T$ act totally on $Y$, determined by $n^x = n+1$, $n \in \ZZ$.  Then $S = W(T,Y)$ is an almost perfect  restriction semigroup. (In fact, by Corollary~\ref{factorizable} below, $S$ is almost factorizable.) Note that each of its $\sigma$-classes is isomorphic to $Y$, as a poset.  Now the monoid $S^1$ is again almost perfect (by part 2 of the cited proposition) but is not perfect, since the only $\sigma$-class with a maximum element is $P_S \cup \{1\}$.\eop\\

Given any semilattice $Y$, a range of examples may be constructed from $T_Y$ itself, for instance by considering $TI_Y$, or any restriction subsemigroups that contain $T_Y$, as `plain' monoids, with $\alpha$ the identical map.  For example, $(T_Y)_{TI_Y, TI_Y}$ is such a semigroup.

Examples of $T$-proper monoids and almost $T$-proper semigroups were given briefly at the end of Section~\ref{properness}.  We shall consider factorizability in the next section.\\

\noindent {\bf Relatively free restriction semigroups}\\

Let $FR_X$ and $FRM_X$ be the free restriction semigroup and monoid respectively, on the set $X$.  It is easily seen (and well known) that $FRM_X = FR_X^1$. The submonoid $T$ of $FRM_X$ generated by $X$ $P$-generates $FRM_X$.   Let $X^*$ be the free monoid on $X$. The natural map $FRM_X \longrightarrow X^*$ restricts to an isomorphism on $T$ and induces $\mathrel{\sigma}$.  Therefore $FRM_X$ is $X^*$-proper and so $FR_X$ is almost $X^*$-proper (and so almost perfect).
We omit further reference to $FR_X$. The semilattice of projections of $FRM_X$ is isomorphic to the semilattice of idempotents of the free {\em inverse\/} monoid on $X$.  While this of course follows from the published structure theorems for $FRM_X$, it can be independently proven (e.g. \cite{jones_P-varieties}).

\begin{corol}\label{free}  The free restriction monoid $FRM_X$ on $X$ is isomorphic to $W(X^*,Y)$, where $X^*$ is the free monoid on $X$, acting on the semilattice of projections $Y$ of $FRM_X$ according to the Munn representation.
\end{corol}

This result can be extended to relatively free restriction semigroups in certain varieties of restriction semigroups.  (See \cite{cornock, jones_restriction1} for background on such varieties.)  Generalizing results on  varieties of inverse semigroups \cite[Section XII.9]{petrich}, Cornock \cite{cornock} showed that a variety $\BV$ of restriction semigroups has {\em proper covers\/} (that is, every member of $\BV$ has a proper cover in $\BV$) if and only if its (relatively) free objects 
 are proper; and that, given any subvariety $\BN$ of the variety $\BM$ of all monoids, regarded as restriction monoids, the class of all restriction semigroups having a proper cover over a member of $\BN$ is a variety having proper covers.

Now by Corollary~\ref{special cover}, any proper restriction semigroup $S$ has an almost $T$-proper cover belonging to the variety generated by $S$.  By standard arguments, the relatively free semigroups in any variety $\BV$ with proper covers must themselves be almost $T$-proper and are therefore representable in the form $W(T,Y)$, where $T$ is the free monoid in $\BV \cap \BM$.  Of course, once again the semilattice of projections needs to be constructed before this representation becomes a concrete structure theorem.


\section{(Almost) factorizable restriction semigroups}\label{embedding theorem}

In any restriction monoid $M$, the $\RP$-class $\RP_1$ is a submonoid, since $\RP$ is a left congruence.  Such a monoid is {\em left factorizable\/} \cite{EQ, GS, S1} if $M = P_M \RP_1$. If $M$ is proper, then $\mathrel{\sigma}$ separates $\RP_1$.  Therefore the proper left factorizable monoids are precisely the $\RP_1$-proper monoids, in our language. 
 Similarly, the proper right factorizable monoids are the $\LP_1$-proper monoids and  
the proper factorizable monoids are the $\HP_1$-monoids.

According to \cite{GS, S1}, a restriction {\em semigroup\/} $S$  is {\em almost left factorizable\/} if 
every element of $ S$ belongs to some member of the $\RP$-class of the identity in $C(S)$.   
Again there are naturally right and two-sided versions of this definition. They extend to restriction semigroups the definition for inverse semigroups of Lawson \cite{lawson_book}.  
The first equivalence in the following result was no doubt known to the authors of \cite{GS, S1}.

\begin{lemma}\label{alf} A proper restriction semigroup $S$ is almost [left, right] factorizable if and only if $C(S)$ is [left, right] factorizable, and thus if and only if $S$ is almost $T$-proper, where $T = $ {\rm[}$\RP_1$, $\LP_1${\rm ]} $\HP_1$ of $C(S)$.
\end{lemma}

{\bf Proof}. Suppose $S$ is almost left factorizable. Let $B \in C(S)$. By hypothesis, each $a \in B$ belongs to some $A \in \,\RP_1$. Since both $B$ and $A$ are contained within $\mathrel{\sigma}$-classes of $S$, in fact $A$ is the same for each choice of $a$, that is $B \subseteq A$.  By Lemma~\ref{basic}, $C(S)$ is $\RP_1$-proper and, by the discussion above, left factorizable. 
Conversely, if $C(S)$ is $\RP_1$-proper, then for each $a \in S$, $a\darrow \subseteq A$, that is, $a \in A$, for some $A \in\, \RP_1$.
The other cases are similar.\eop\\

The semigroups $W(T,Y)$ that arise in the representations of the left and two-sided almost factorizable semigroups in the following result are precisely the `original' ones of  \cite{FG, GS, S1} (see the discussion following Theorem~\ref{X semigroup}).  Note that the right-hand version falls within the realm of the broader $W$-product, but not within that of the original.  The rather simpler statement of the proposition in the case of factorizability is left to the reader.

\begin{propo}\label{factorizable} The following are equivalent for a proper restriction semigroup $S$:

 \be[(a)]
 \item $S$ is almost left factorizable {\rm [}almost right factorizable, almost factorizable{\rm ]};
 
 \item $S$ is almost $T$-proper, where $T$ is the $\RP$-class {\rm [}$\LP$-class, $\HP$-class{\rm ]} of the identity in $C(S)$;
  
 \item the action that is induced by the Munn representation of $S$, by isomorphisms between ideals of $P_S$,  is by endomorphisms {\rm [}onto mappings, automorphisms{\rm ]};
  
  \item $M \cong W(T,Y)$ for such an action of a monoid $T$ upon a semilattice $Y$ (that is, in the left and two-sided cases, the `original' $W$-product).
  
 
  \ee
\end{propo}

{\bf Proof.} First consider the left-hand case.  The first equivalence was shown in Lemma~\ref{alf}.  Now suppose $C(S)$ is $\RP_1$-proper.  Then the action induced by the Munn representation of $S$,  according to Proposition~\ref{C(S)} and its corollary, is by fully defined mappings, that is, by endomorphisms.  An application of Theorem~\ref{converse} yields (d).  That $W(T,Y)$ as in (d) is almost left factorizable can be checked directly, or by using Proposition~\ref{X action} to prove (b).

In the right-hand case, by duality, the action is instead by {\em onto\/} partial mappings.  In the two-sided case, the requirement becomes that the action be by {\em automorphisms\/} of $Y$.\eop\\

That any restriction semigroup with one of the factorizability properties considered above has a proper cover of the same type was shown in \cite{GS, S1}. This can also be obtained by direct calculation using Corollary~\ref{T-proper cover}, the simplest form of the coverings in that section.\\

Szendrei \cite[Theorem 4.1]{S2} proved that any restriction semigroup is embeddable in an almost left factorizable semigroup.   In a paper that the author received as this work was nearing completion, 
Kudryavtseva \cite{K} provided an alternative proof as a consequence of more general embedding results (see the last paragraphs of Section~\ref{intro}).  Gomes, Gould and Hartmann (personal communication) have proved that in fact any restriction semigroup is embeddable in an almost factorizable semigroup (and therefore, by definition, into a factorizable monoid).  In this section, we apply our methods to provide a quite different proof of their result, based solely on the following result of McAlister.  

\begin{result}\cite[Corollary 4.4]{McA}\label{mcalister} Any semilattice $Y$ is embeddable as an ideal in a semilattice $Z$ with the property that any isomorphism between principal ideals of $Z$ can be extended to an automorphism of $Z$.
\end{result}

After seeing the preprint \cite{K} by Kudryavtseva, we split the proof into two parts in order to illustrate how her Theorem 28 could also be recovered, in an apparently slightly stronger form.

\begin{propo}\label{kudry}(cf \cite[Theorem 28]{K})  Let $Y$ be a semilattice with identity and $X^*$ a free monoid that acts on $Y$ by isomorphisms between principal ideals.  
Then $Y$ is embeddable as an ideal in a semilattice $Z$, with identity, and the action of $X^*$ on $Y$ extends to an action on $Z$ by automorphisms, in such a way that $W(X^*, Y) $ is embeddable in the `original' $W$-product $W(X^*, Z)$ (in the sense of \cite{FG}).  The latter semigroup is a factorizable restriction monoid.
\end{propo}

{\bf Proof}.  Let $Z$ be a semilattice as in Result~\ref{mcalister}; it may be assumed that $Z$ has a greatest element (by adjoining one if that is not the case) and that $Y$ itself is an ideal of $Z$.  Then any principal ideal of $Y$ is also a principal ideal of $Z$ and so $T_Y$ is an inverse subsemigroup of $T_Z$.

Put $T = X^*$.  For each $\delta \in T_Y$, denote by $\overline{\delta}$ a specific extension to an automorphism of $Z$. In particular, for each $x \in X$, $ \overline{x\alpha} \in {\rm Aut} Z$.  The map $ x \mapsto \overline{x\alpha}$ extends uniquely to a monoidal homomorphism $\oalpha: T \longrightarrow {\rm Aut} Z$.

By Theorem~\ref{X semigroup}, $W(T, Z)$ is a perfect restriction monoid.  By  Proposition~\ref{factorizable}, it is factorizable.  If $e \in Z$ and $t \in T$, then by a slight abuse of notation we shall write $e^{\overline{t}}$ for $e (t \oalpha)$.
Suppose $t = x_1 \cdots x_n$ and $e \in \Delta t \subseteq Y$. Then $e \in \Delta x_1$ and so $e (\overline{x_1 \alpha}) = e (x_1 \alpha) \in \Delta x_2$.  An inductive argument yields that $e (t \oalpha) = e (t \alpha)$ or, in the above notation, $e^{\overline{t}} = e^t$, that is, $t\oalpha$ extends $t \alpha$ for all $t \in T$.

It follows that $W(T,Y)$ is a restriction subsemigroup of $W(T, Z)$ (but not in general a submonoid).\eop

\begin{theorem}\cite{gould et al}\label{embedding} Any restriction semigroup is embeddable in a factorizable restriction monoid.
\end{theorem}

{\bf Proof}.  The outline of the proof is as follows.  Given the restriction semigroup $R$, construct an almost perfect cover $W(T,Y)$ as in Corollary~\ref{T-proper cover alt}. By the previous proposition, this cover embeds in $W(T, Z)$, where $Z$ is as in Result~\ref{mcalister}, and this semigroup is factorizable.  The congruence induced on $W(T,Y)$ by the covering map is explicitly defined and has a natural extension to $W(T,Z)$.  Thus $R$ embeds in a quotient of the latter, which is again factorizable. Here are the details.

By adjoining an identity, if necessary, it may be assumed that $R$ is a monoid.  Let $X$ generate $R$, as such, let $T = X^*$ and $\phi: T \longrightarrow R$ the natural homomorphism.  Construct the covering $R_T$ according to Corollary~\ref{T-proper cover}, interpreted as in Corollary~\ref{T-proper cover alt}.  

That is, $R_T \cong W(T,Y)$, where $Y = P_R$ and the action is induced by the composition $\alpha = \phi \theta: T \longrightarrow T_Y$, $\theta$ being the Munn representation of $R$ in $T_Y$, that is, for $t \in T$ and $e \in \Delta t$, $e^t = (e(t\phi))^*$.  By the remarks following the latter corollary, the covering map $W(T,Y) \longrightarrow R$ is given by $(t, f) \mapsto (t\phi) f$ or, equivalently, by $(t, e^t) \mapsto e (t\phi)$.

We now use the extension of the action specified in the proof of Proposition~\ref{kudry}, without further comment.

Define a relation $\sim$ on $W(T,Z)$ by 
\[(t, f) \sim (u,h) \text{  if  } t\phi((t\phi)^* f) = u\phi( (u \phi)^* h) \text{ in } R.\]
Although $f, h \in Z$, $ (t \phi)^*f , (u \phi)^* h \in Y$, since $Y$ is an ideal of $Z$.  If $(t,f), (u,h) \in W(T,Y)$, then $(t,f ) \sim (u,h)$ if and only if $(t \phi) f = (u \phi) h$, that is, the pairs are identified under the covering map $W(T,Y) \longrightarrow R$.  Therefore, once we have proved that $\sim$ is a congruence on $W(T,Z)$, the quotient monoid $W(T,Z)/\sim$ will contain $ R$.

We first prove (cf the alternative formula for the covering map) that $t \phi ((t\phi)^*e^{\overline{t}}) =  (e (t\phi)^+) t\phi$ for any $(t, e^{\overline{t}}) \in W(T,Z)$. 
Now, recalling the definition of the action of $T$ on $Y$ stated above,  \[(t\phi)^*e^{\overline{t}} = ( (t \phi)^+)^t e^{\overline{t}} =  ( (t \phi)^+)^{\overline{t}} e^{\overline{t}} =  ( (t \phi)^+ e)^{\overline{t}} =  ( (t \phi)^+ e)^t = ( ( (t \phi)^+ e) (t\phi))^*.\]  
Therefore $t \phi ((t\phi)^*e^{\overline{t}})  = (t \phi)  ( ( (t \phi)^+ e) (t\phi))^* = ( (t \phi)^+ e) t \phi$, using the identity $y(xy)^* = x^*y$ (the dual of the last of the four identities presented in Section~\ref{prelim}).

Since any element $(s,k) $ of $W(T,Z)$ can be factored as $(1,k)(s,1)$, it suffices to show that $\sim$ is respected by multiplication on the left and right by $(s,1)$ and $ (1,k)$ separately. We will verify this on the left, the other case being similar (indeed, in view of the formula proved above, dual).

Suppose $(t,f) \sim (u,h)$, so that $t \phi( (t \phi)^* f) = u \phi ( (u \phi)^* h)$. 
First consider left multiplication by $(s,1)$.  Here $(s,1)(t, f) = (st, f)$ and $(s,1)(u, h) = (su, h)$.  
Now since $(st)\phi = (s\phi)(t\phi)$, $(st)\phi = (st)\phi (t\phi)^*$, so
\[(st)\phi ( ((st) \phi)^* f) =  (st)\phi ((st) \phi)^*((t \phi)^*f) = (s\phi)(t\phi)((t \phi)^*f).\]  
Similarly, $(su)\phi (((su)\phi)^* h) = (s\phi)(u\phi) ((u\phi)^*  h)$. 
Substituting from the equation for $(t,f) \sim (u,h)$ above yields that $(st,f) \sim (su, h)$.

Next consider left multiplication by $(1,k)$ and use the `dual' form,  that is, write $(t,f) = (t, e^{\overline{t}})$ and $(u,h) = (u, g^{\overline{u}})$.  Using the alternative formula above,  
 $(1,k) (u, g^{\overline{u}}) =  (u, (kg)^{\overline{u}})$, likewise.    Then  \[(ke (t\phi)^+) t \phi =   (k (u\phi)^+) (e (t\phi)^+) t \phi= (k (u \phi)^+) (g (u\phi)^+) u \phi = (kg (u\phi)^+) u\phi\] and therefore $ (t, (ke)^{\overline{t}}) \sim  (u, (kg)^{\overline{u}})$.

The proof is completed on noting that a monoidal quotient of a factorizable restriction monoid is again factorizable (cf \cite[Lemmas 4.7, 5.8]{GS}).\eop

\section{Proper restriction semigroups in general}\label{proper}

As noted in Section~\ref{intro}, only minor modifications, based on Lemma~\ref{sub} and the general case of Corollary~\ref{T-rep}, are needed to extend all of the main results of this paper to proper restriction semigroups (and monoids).  In this section, we focus on the extensions of each of the cited results to the proper case.  We refer the reader to the statement of Theorem~\ref{newcovers} and the construction that precedes it.  We should note here that, in the general situation of that theorem, the connection between $S$ and $R$ may be too tenuous to deduce the converse statement in the next result.  Instead, therefore, we phrase it in terms of Corollary~\ref{C(S) cover}.  As noted there, in view of Theorem~\ref{P-sep covers} that situation includes the case $R = S^1$ and, moreover, is in a sense general. 

\begin{theorem} \label{proper newcovers} (cf Theorem~\ref{newcovers}) If $\alpha: T \longrightarrow R$ is a subhomomorphism, then $S_{T,R}$ is a proper restriction semigroup.  

If $S \leq R \leq C(S)$ (cf Corollary~\ref{C(S) cover}), then $S_{T,R}$ is almost perfect if and only if $\alpha$ is a homomorphism.  In particular, this is true for $S_T$.

If $S$ is a monoid and $\alpha$ is monoidal, then $S_T$ is an $F$-restriction monoid;
thus $S_T$ is $T$-proper (that is, perfect) if and only if $\alpha$ is also a homomorphism.

\end{theorem}

{\bf Outline of proof}.  It was mentioned during the proof of the cited theorem that only the subhomomorphism property of $\alpha$ was required for closure; the homomorphism property was not used in the proof of properness.  In the monoidal case, the stated property of $\sigma$-classes was also proved there.

Next we show that in the case $S \leq R \leq C(S)$, if $\sigma$ is perfect then $\alpha$ must be a homomorphism.  For clarity's sake, we shall explicitly use the embedding $\tau_S $, where $a \tau_S = a\darrow$.   In the notation of the cited theorem, the $\sigma$-classes of $S_{T,R}$ are the sets $S_t$, $t \in T$, where $S_t = \{(a \darrow ,t) : a \darrow\leq t\alpha\} \text{ in } C(S)$.  Note that if $A \in C(S)$ and $a \in S$, then $a \darrow \leq A$ if and only if $a \in A$.

Let $t, u \in T$.  Each of $t\alpha$, $u\alpha$ and $(tu)\alpha$ is an order ideal of $S$. 
It must be shown that 
 $(t u)\alpha \subseteq (t \alpha)(u\alpha)$ in $C(S)$.  Let $a \in (tu) \alpha$.  So $a \darrow \leq (tu)\alpha$ and therefore $(a \tau_S, tu) \in S_{tu}$.  By perfection of $\sigma$, $(a \tau_S,tu) = (b \tau_S,t) (c \tau_S,u)$ for some $b \in t\alpha$, $c \in u\alpha$.  Since $\tau_S$ is injective, $a = bc \in (t\alpha)(u\alpha)$, as required.
 
 In the monoidal case, either the proof above may be modified or one can directly use the assumption that the  elements $(t\alpha, t)$, $t \in T$ form a subsemigroup to show that $\alpha$ is a homomorphism.\eop\\

With `almost perfect' replaced by `proper' in the semigroup case, and `perfect' replaced by `$F$-restriction'
in the monoid case (so that the map $\kappa: T \longrightarrow N$ is defined), Theorem~\ref{P-sep covers} and Corollary~\ref{all covers} now characterize the respective covers in terms of the broader construction in the theorem just stated.

In the context of this section, if there is a subhomomorphism $\alpha: T \longrightarrow TI_Y$, then $T$ is said to {\em sub-act\/} on the semilattice $Y$.  We otherwise retain the notation of Section~\ref{X}.  The construction $W(T, Y)$ is as in that section; as noted at that point, the argument for closure only required a sub-action.

\begin{theorem}\label{proper X semigroup} (cf Theorem~\ref{X semigroup}) Let $T$ be a monoid, $Y$ a semilattice, and $\alpha : T \longrightarrow TI_Y$ a monoidal subhomomorphism, that is, $T$ sub-acts on $Y$ by isomorphisms between ideals.  Then $W = W(T,Y)$ is isomorphic to the semigroup $(T_Y)_{T, TI_Y}$ and so 
is a proper restriction semigroup, with $P _W \cong Y$ and $W /\sigma \cong T $; it is almost $T$-proper if and only if $\alpha$ is a homomorphism.

Further
 if $Y$ is a monoid and $\alpha$ is a subhomomorphism into $T_Y$, that is, $T$ sub-acts on $Y$ by isomorphisms between principal ideals, then $W$ is isomorphic to the monoid $(T_Y)_T$ and so is an $F$-restriction monoid;
 it is $T$-proper if and only if $\alpha$ is a homomorphism.
\end{theorem}

{\bf Proof}.  This follows from Theorem~\ref{proper newcovers} and the cited theorem.\eop\\

Finally, the converse again follows from application of the above and the cited theorem.

\begin{theorem}\label{proper converse} (cf Theorem~\ref{converse}) Let $S$ be a proper restriction semigroup.  Put $T = S/\sigma$ and $Y = P_S$.  Then $S \cong W(T,Y)$, where the sub-action of $T$ on $Y$, by isomorphisms between ideals, is induced by the subhomomorphism $\kappa \otheta$, the composition of the injection of $T$ in $C(S)$ with the extension of the Munn representation of $S$ to $C(S)$.

If $S$ is an $F$-restriction monoid
then $Y = Y^1$ and the sub-action of $T$ on $Y$ is by isomorphisms between principal ideals, induced by the Munn representation of $S$ itself.
\end{theorem}

Cornock and Gould \cite{CG} provided a structure theorem for proper restriction semigroups in general, based on {\em pairs\/} of actions of a monoid $T$ on a semilattice $Y$.  For such a pair, they define a semigroup ${\cal M}(T, Y)$ and show that this construction describes proper restriction semigroups.  Clearly, there must be a correspondence between their construction and that in this section, but we have not pursued this explicitly since the general case is not the one of main interest in our work.
\section{Specialization to inverse semigroups}\label{inverse}
Recall from Section~\ref{prelim} that proper inverse semigroups are usually termed $E$-unitary. The specializations of the general parts of Theorems~\ref{proper X semigroup} and \ref{proper converse} to inverse semigroups are easily obtained.  We use the terminology of \cite{petrich}, that a {\em prehomomorphism\/} of inverse semigroups is a subhomomorphism that respects inverses.  We leave the statement in the case of monoids to the reader.

\begin{corol}\label{inverse converse} (cf \cite{PR}) Let $T$ be a group, $Y$ a semilattice, and $\alpha : T \longrightarrow TI_Y$ a monoidal prehomomorphism.  Then $W = W(T,Y)$ is an $E$-unitary inverse semigroup, 
isomorphic to $(T_Y)_{T, TI_Y}$, with 
$E _W \cong Y$ and $W /\sigma \cong T $.  Conversely, for any $E$-unitary inverse semigroup $S$, let $T$ be the group $S/\sigma$ and $Y = E_S$.  Then $S \cong W(T,Y)$, where the sub-action of $T$ on $Y$, by isomorphisms between ideals, is induced by the prehomomorphism $\kappa \otheta$.
\end{corol}

{\bf Proof}. In the direct case, the conclusion is clearer if we consider the semigroups $S_{T,R}$ of Theorem~\ref{X semigroup}, with $R$ and $S$ inverse semigroups, under the same assumptions on $T$ and $\alpha$.  (We may then quote the second isomorphism in the theorem.)  If $(s,t) \in S_{T,R}$, then $s \leq t \alpha$; by the assumption on $\alpha$ and compatibility of the natural partial with inverses in  inverse semigroups, $s^{-1} \leq t^{-1} \alpha$, so $(s^{-1}, t^{-1})$ is an inverse for $(s,t)$.

In the converse case, $T$ is necessarily a group and $P_S = E_S$.  Since $(a\sigma)^{-1} = a^{-1}\sigma$ in any inverse semigroup, $\kappa$, and therefore $\kappa\otheta$, is a prehomomorphism.\eop\\

The theory of $E$-unitary inverse semigroups is very well mined and it is hardly to be expected that the general theory, such as this corollary, would reduce to anything but well-trodden ground (albeit in somewhat different language), even though that ground played no part in the development of our theory.  As cited incidentally above,  the main body of Corollary~\ref{inverse converse} above is essentially the description of $E$-unitary inverse semigroups found by Petrich and Reilly \cite{PR} (see also \cite[Theorem VI.8.12]{petrich}, where the inverse semigroup of isomorphisms between ideals of a semilattice $Y$ is denoted $\Sigma(Y)$).

Since we are more interested in where our theory does not simply extend inverse semigroup theory, we refer the reader to Chapter VII of the monograph \cite{petrich} by Petrich for comprehensive coverage. Of course, McAlister's $P$-theory is the gold standard, especially when treated in concert with the various alternative descriptions of $E$-unitary covers discovered in subsequent years (\cite[Section VII.4]{petrich}), where the semigroup $C(S)$ not surprisingly plays a significant role.  Chapter 7 of the monograph \cite{lawson_book} by Lawson presents the $P$-theorem from a somewhat different perspective.

The point that we wish to emphasize here is that the specialization of almost perfection and perfection to inverse semigroups does not yield a general theory as it does for restriction semigroups.  They are well-studied classes:

 \begin{propo} \label{specialization}  The almost perfect inverse semigroups are the semidirect products of semilattices and groups. The perfect inverse monoids are the monoidal such products. 
  \end{propo}
 {\bf Proof}.  This is simply the combination of Propositions~VII.5.11, VII.5.14, and VII.5.24 of \cite{petrich}.\eop\\
 
 The key point of divergence in our work is then the covering result Theorem~\ref{newcovers} and, more particularly, the simpler result Corollary~\ref{T-proper cover}.  The latter produces an almost perfect cover $S_T$ for $S$ which, as we have just seen, is a semidirect product of a semilattice and a group.  But \cite[Theorem 7.10]{lawson_book} only the almost factorizable inverse semigroups are quotients of such semidirect products.  In the monoidal case, only the factorizable inverse semigroups are quotients.


\noindent
Department of Mathematics, Statistics and Computer Science\\
Marquette University\\
Milwaukee, WI 53201, USA\\
peter.jones@mu.edu

\begin{thebibliography}{99}
\bibitem{BGG}
	M. Branco, G. Gomes and V. Gould, Left adequate and left Ehresmann monoids, Internat. J. Algebra Comput. 	21 (2011), 1259--1284.
	
\bibitem{cornock}
	C. Cornock, Restriction Semigroups: Structure, Varieties and Presentations, Ph.D. Thesis, York, England, 	2011.
	
\bibitem{CG}
	C. Cornock and V. Gould, Proper two-sided restriction semigroups and partial actions, J. Pure Appl. Algebra 	216 (2012), 935--949.

\bibitem{EQ}
	A. El Qallali, Factorisable right adequate semigroups, Proc. Edinburgh Math. Soc. 24 (1981), 171--178.
	
\bibitem{FG}
	J. Fountain and G. Gomes, Proper left type-A monoids revisited, Glasgow Math. J. 35 (1993), 293--306.
	
\bibitem{FGG1}
	J. Fountain, G. Gomes and V. Gould, A Munn type representation for a class of $E$-semiadequate semigroups, J. Algebra 218 (1999), 69 --714.
	
\bibitem{FGG}
	J. Fountain, G. Gomes and V. Gould, The free ample monoid, Internat. J. Alg. Comp. 19 (2009), 527--554.
	
\bibitem{GG}
	G. Gomes and V. Gould, Fundamental Ehresmann semigroups, Semigroup Forum 63 (2001), 11--33.
	
\bibitem{GS}
	G. Gomes and M. Szendrei, Almost factorizable weakly ample semigroups, Comm. Algebra 35 (2007), 3503--3523.
	
\bibitem{GG2}
	G. Gomes and V. Gould, Left adequate and left Ehresmann monoids II, J. Algebra 348 (2011), 171--195.	
	
\bibitem{gould et al}
	G. Gomes, V. Gould and M. Hartmann, personal communication.
	
\bibitem{howie}
	John M. Howie, Fundamentals of Semigroup Theory, Clarendon Press, Oxford, 1995.
	
\bibitem{jones_P-restriction}
	Peter R. Jones,  A common framework for restriction semigroups and regular $*$-semigroups, J. Pure Applied Algebra 216 (2012), 618--632.
	
\bibitem{jones_restriction1}
	Peter R. Jones, On lattices of varieties of restriction semigroups, Semigroup Forum 86  (2013), 337--361.
	
\bibitem{jones_P-varieties}
	Peter R. Jones, Varieties of $P$-restriction semigroups, Comm. Algebra 42 (2014), 1811--1834.
	
\bibitem{K}
	G. Kudryavtseva, Partial monoid actions and a class of restriction semigroups, preprint.
	
\bibitem{lawson}
	M. Lawson, Semigroups and ordered categories I: the reduced case, J. Algebra 141 (1991), 422--462.
	
\bibitem{lawson_book}
	M. Lawson, Inverse Semigroups: The Theory of Partial Symmetries.  World Scientific, Singapore, 1998.
	
\bibitem{McA}
	D.B. McAlister, $E$-unitary inverse semigroups over semilattices, Glasgow Math. J. 19 (1978), 1 -- 12.


	

	
\bibitem{petrich}
	M. Petrich, Inverse semigroups, Wiley, New York, 1984.
	
\bibitem{PR}
	M. Petrich and N. Reilly, A representation of $E$-unitary inverse semigroups, Quart. J. Math. Oxford 30 			(1979), 339--350.
	
\bibitem{S1}
	M. Szendrei, Embedding into almost left factorizable restriction semigroups, Comm. Algebra 41 (2013), 1458--1483.
	
\bibitem{S2}
	M. Szendrei, Proper covers of restriction semigroups and W-products, Internat. J. Algebra Comput. 22 (2012), 


\end{thebibliography}
\end{document}